\documentstyle[12pt]{article}
\title{Metric Entropy of Homogeneous Spaces and Finsler Geometry of
Classical Lie Groups
\thanks{Partially supported by a grant from the
National Science Foundation.}}
\author{Stanislaw J. Szarek \\ Department of Mathematics \\
Case Western Reserve University \\ Cleveland Ohio 44106}
\date{}

\begin{document}
\maketitle

\begin{abstract}
\noindent For a (compact) subset $K$ of a metric space and $\varepsilon >
0$, the
{\em covering number}
$N(K , \varepsilon )$ is defined as the smallest number of balls of
radius $\varepsilon$ whose union covers $K$.  Knowledge of the {\em metric
entropy},
i.e., the asymptotic behaviour of
covering numbers for (families of) metric spaces is important in many areas of
mathematics (geometry, functional analysis, probability, coding theory, to
name a few).
In this paper we give asymptotically correct estimates for covering numbers
for a large class of
homogeneous spaces of unitary (or orthogonal) groups with respect to some
natural metrics,
most notably the one induced by the operator norm.
This generalizes earlier author's results concerning covering numbers of
Grassmann manifolds; the generalization is motivated by applications to
noncommutative probability and operator algebras. In the process we give a
characterization
of geodesics in $U(n)$ (or $SO(m)$) for a class of non-Riemannian metric
structures.
\end{abstract}

\newpage
\section{INTRODUCTION}

If $(M, \rho )$ is a metric space,  $K \subset M$ a compact subset and
$\varepsilon > 0$, the
{\em covering number}
$N(K, \varepsilon ) = N (K, \rho , \varepsilon )$ is defined as the
smallest number of balls of
radius $\varepsilon$ whose union covers $K$.   If $K$ is a ball of radius
$R$ in a normed space
of real dimension $d$,  it is easily shown (by a volume comparison
argument) that,
for any $\varepsilon \in (0, R]$,

\begin{equation}
(R / \varepsilon )^d \leq \ N(K, \varepsilon ) \leq (1 + 2R/\varepsilon )^d.
\end{equation}

The lower and the upper estimate in (1) differ roughly by a factor of
$2^d$,  and for
many applications such an accuracy is sufficient.  On the other hand,
determining more
precise asymptotics for covering numbers and their "cousins",  {\em packing
numbers}  (see section 2),
e.g. for Euclidean balls is a
nontrivial proposition (and a major industry).  In this paper we attempt to
obtain
estimates of type (1) for homogenous spaces of the
orthogonal group $SO(n)$ or the unitary group $U(n)$  (like e.g. the
Grassmann manifold
$G_{n,k}$ or the Stieffel manifold, equipped with some natural
metric; we admit metrics induced by unitary ideal norms of matrices,  most notably the
operator norm).  A typical result will be:  if $M$ is a ``nice" homogenous
space of
$SO(n)$ or  $U(n)$ and $\varepsilon \in (0, \theta(M)]$  (where $\theta(M)$
is some
{\em computable} ``characteristic" of $M$,  in the more ``regular" cases
$\theta(M) \approx diam \, M$),  then

\begin{equation}
(c \, diam \, M/ \varepsilon )^d \leq N(M, \varepsilon) \leq (C \, diam\,
M/\varepsilon )^d,
\end{equation}

\noindent where $d$ is the (real) dimension of $M$  and  $c$ and $C$ are
constants
independent of
$\varepsilon $ and  (largely) of $M$.
Of course universality of the constants in question is {\em the} crucial point.

We note in passing that (again,  by a standard ``volume comparison"
argument) (2) is
equivalent to the (normalized)
Haar measure of a ball of radius $\varepsilon$ being between
$(c_1 \varepsilon /diam \, M)^{\dim M}$ and $(C_1 \varepsilon /diam \,
M)^{\dim M}$,
where $c_1, C_1 > 0$ are {\em universal} constants related to $c, C$.  We
also point
out that the underlying metric being typically non-Riemannian,  the methods of
Riemannian geometry do not directly apply.

This paper is an elaboration of the note  \cite{Iowa} by the author,  where
(2) was proved
in the special case of $M$
being $SO(n)$,   $U(n)$ or a Grassmann manifold $G_{n,k}.$  (Another
argument was
given later in \cite{pajor}.) The original motivation and application was the
finite-dimensional basis problem;  more precisely, (2) was used in the proof of
(Theorem 1.1 in \cite{fdbp}):

\medskip\noindent {\em There is a constant} $c > 0$
{\em such that,
for every positive integer }$n$, {\em there exists an } $n${\em-dimensional
normed space} $B$
{\em such that, for every projection} $P$ {\em on }$B$ {\em with, say, }
$.01 n \leq rank P \leq .99 n${\em , we have }
$\| P:B \rightarrow B \| > c \sqrt{n}$.

\medskip The results of \cite{Iowa} were subsequently applied to other
problems in convexity,
local theory of Banach spaces,  operator theory,  noncommutative probability
and operator algebras
(cf. e.g. \cite{herrero}, \cite{exotic}, \cite{cube}, \cite{voiculescu}).
It turned out recently (see \cite{voiculescu}) that some questions from the
last
two fields lead naturally to queries about validity of
estimates of type (2)  in settings more general than that of \cite{Iowa}
(which,
additionally,  had a rather limited circulation).  It is the purpose of
this paper to
provide a {\em reasonably general} answer to such questions.   However,
describing
asymptotics for {\em completely general}  homogenous spaces of $SO(m)$ or
$U(n)$ is,
in all likelihood,  hopeless.  In this paper we cover a number of special
cases,
including those that
have been explicitly inquired about.   We identify (easily computable)
invariants
relevant to the problem and provide ``tricks" that could be potentially useful
to handle cases not addressed here.  Ecxept for brief comments here and there,
we restrict our attention to $M = G/H$,
where $G = SO(m)$ or $U(n)$, and $H$ - a connected Lie subgroup of $G$;
but clearly
most of our  analysis can be extended to other compact linear Lie groups.

The organization of the paper is as follows.

Section 2 explains notation and presents various preliminary results concerning
covering numbers and their relatives,  unitarily invariant norms and the
exponential map.
Most of these are known or probably known.  In section 3 we show
(in Proposition 6) that cosets of
one parameter semigroups are the geodesics in $U(n)$ (or $SO(m)$) endowed
with an
intrinsic metric induced by a unitarily invariant norm;  this result
could be of independent interest. In section 4 we discuss several simple
examples
that point out possible obstructions to estimates of type (2) and suggest
invariants
mentioned above. We then use the results of the preceding two sections to show
estimates of type (2) for an abstract class of homogeneous spaces
that contains $U(n)$, $SO(n)$ and $G_{n,k}$ (Theorem 8).  In section 5 we
discuss various
possibilities for relaxation of assumptions from section 4,  in particular we
cover the special cases motivated by applications.  At the end of section 5
we briefly address the issue of extending the results to metrics generated
by unitarily invariant norms other than the operator norm.

{\em Acknowledgement.} The author would like to  express his gratitude to
D. Voiculescu,  whose
encouragement was instrumental both in the inception and the completion
of this work.  The final part of the research has been performed while the
author was in
residence at MSRI Berkeley; thanks are due to the staff of the institute
and the organizers of the Convex Geometry semester for their hospitality
and support.

\section{NOTATION AND PRELIMINARIES}

We start with several remarks clarifying the relationship between covering
numbers,
packing numbers and their slightly different versions that exist on the
market.  First, since
the centers of balls in the definition of $N(K,\varepsilon)$, as given in
the introduction, do not
necessarily need to be in $K$,  the exact value of $N(K, \cdot)$ may depend
on the ambient
metric space $(M,\rho)$ containing $K$. Accordingly,  it is sometimes more
convenient to allow sets of
{\em diameter} $\leq 2 \varepsilon$ in place of balls
of {\em radius} $\varepsilon$;  call the resulting the quantity
$N'(K,\varepsilon)$.  If the
centers of the balls in the definition {\em are} required to be in $K$,
call the quantity
$N''(K, \varepsilon)$.  Finally,  let the {\em packing number}
$\tilde{N}(K,\varepsilon)$ be
defined as the maximal cardinality of an
$\varepsilon$-separated  (i.e.  $\rho(x,x') > \varepsilon$ if $x \neq x'$)
set in $K$.  The quantities $N$, $N'$ and $\tilde{N}$ are related as follows

\begin{equation}
N'(\cdot,\varepsilon) \leq N(\cdot,\varepsilon) \leq N''(\cdot,\varepsilon)
\leq \tilde{N}(\cdot,\varepsilon)
\leq N'(\cdot,\varepsilon/2).
\end{equation}

\noindent Consequently,  for our ``asymptotic" results the four quantities
are essentially interchangeable.

 If the metric space  $(M,\rho)$ is
actually a normed space with a norm $\| \cdot \|$ and unit ball  $B$,
we may write $N(K, \| \cdot \| , \varepsilon )$ or $N(K, B , \varepsilon )$
instead of $N(K, \varepsilon )$ or $N(K, \rho , \varepsilon )$. The technology
for estimating covering/packing
numbers of subsets (particularly convex subsets) of normed spaces is quite
well-developed
and frequently rather sophisticated  (see \cite{carl}, \cite{pisier}).  We
quote here a
simple well-known result
that expresses $N(\cdot , \cdot )$ in terms of a ``volume ratio", and of which
(1) is a special case.

\medskip\noindent {\bf Lemma 1.}  {\em Let } $K$, $B \subset R^d$ {\em with
} $B$
{\em convex symmetric.  Then, for any } $\varepsilon >0$,
$$( \frac{1}{\varepsilon} )^d \frac{vol \, K}{vol \, B} \leq N(K,B ,
\varepsilon )
\leq N''(K,B , \varepsilon ) \leq N'(K,B,\varepsilon/2)
\leq (\frac{2}{\varepsilon} )^d \frac{vol \,(K + \varepsilon/2 B)}{vol \, B}.$$

The next lemma is just an observation which expresses the fact that the
covering/packing
numbers are invariants of Lipschitz maps.

\medskip\noindent {\bf Lemma 2:} {\em Let } $(M, \rho)$ {\em and } $(M_1,
\rho_1)$
{\em be metric spaces, } $K \subset M$, $\Phi : K \rightarrow M_1$, {\em
and let }
$L > 0$.  {\em If } $\Phi$ {\em verifies}
$$\rho_1( \Phi (x), \Phi (y)) \leq L \, \rho (x,y) \, \mbox{ {\em for} } \,
x,y \in K$$
{\em (i.e. } $\Phi$ {\em is a Lipschitz map with constant } $L$),
{\em  then, for every } $\varepsilon > 0$,
$$N''(\Phi (K), \rho_1 , L \varepsilon) \leq N''(K, \rho , \varepsilon ).$$
{\em Moreover, } $N''$ {\em can be replaced by } $N'$ {\em or } $\tilde{N}$
{\em and,  if } $\Phi$ {\em can be extended to a function on }
$M$ {\em that is still Lipschitz with constant } $L$, {\em also by } $N$.

\medskip We now turn to our main interest,  the unitary group $U(n)$, the
(special) orthogonal
group $SO(n)$ and their homogeneous spaces.  (As $O(n)$ is geometrically a
disjoint
union of two copies of $SO(n)$,  all statements about $SO(n)$ will easily
transfer to $O(n)$.)
Throughout the paper we will reserve the letter $G$ to denote,  depending on the context,
 $U(n)$ or $SO(n)$. Similarly,  we will reserve the letter $\cal{G}$ to
denote the
Lie algebra of $G$,  the space $u(n)$ or $so(n)$ of skew-symmetric matrices.
Since $G$ and $\cal{G}$ are subsets of $M(n)$  (the algebra of $n \times n$
matrices,
real or complex as appropriate),  they inherit various metric structures
from the
latter.  In this paper we focus on the one induced by the operator norm
(as an operator on the Euclidean space,  that is),  but will also consider
the Schatten
$C_p$-norms $\|x\|_p = (tr \,|x|^p)^{1/p}$ with the operator norm
$\| \cdot \|_{op} = \| \cdot \| _{\infty}$ being the limit case.  We will
use the same notation
$\|x\|_p$ for the $\ell_p^n$-norm on ${\bf R}^n$ or ${\bf C}^n$,  but this
should not
lead to confusion.  We will also occasionally mention other unitarily
invariant norms
(i.e. verifying $\|x\| = \| uxv \|$ if $x \in M(n)$ and $u,v \in G$), each
necessarily associated
with a symmetric norm on ${\bf R}^n$ (which we will also denote by $\|
\cdot \|$) via
$\|x\| = \|(s_k (x))_{k=1}^n \|$,  where $s_1 (x), \dots, s_n (x)$ are
singular numbers
of $x$.

If $\rho$ is a metric on $G$ and
$H \subset G$ a closed subgroup,  we consider
the homogeneous space $M = G/H$  of left cosets of $H$ in $G$ as endowed with
the canonical quotient metric  $\rho_M (E,F) = \inf \{ \rho (u,v) \, : \, u
\in E, v \in F \}$.
A fundamental example is that of the Grassmann manifold $G_{n,k}$ of
$k$-dimensional subspaces
of ${\bf R}^n$ (resp. ${\bf C}^n)$: the relevant subgroup $H$ of $SO(n)$
(resp. $U(n)$) consists of matrices of the form
\begin{equation}
\left [
\begin{array}{cc}
u_1 & 0 \nonumber \\
0 & u_2
\end{array}
\right ]
\end{equation}
where $u_1 \in SO(k)$ and $u_2 \in SO(n-k)$  (resp. $ U(k), U(n-k)$) and
the identfication
of cosets of $H$ with the subspaces is via $uH \sim uE_k$,  where $E_k$ is
the linear span of the
first $k$ vectors of the standard basis of ${\bf R}^n$ (resp. ${\bf C}^n)$.

It will be usually convenient to consider the {\em intrinsic}
metric,  which we will again call $\rho$, on $M=U(n)$ (or $SO(n)$, or the
homogeneous space):
 $\rho(u, v)$ is the infimum  of
lengths of curves in $M$ connecting $u$ and $v$.  Since this is going to be
relevant later,
we observe that the infimum may be taken over {\em absolutely continuous}
curves
$\gamma \, : \, [a,b] \rightarrow M$ (the infimum is then achieved, while any
rectifiable curve parametrized by arc length is absolutely continuous).
Then the length
\begin{equation}
\ell(\gamma) = \int_0^1 \| \gamma'(t)\|dt,
\end{equation}
where  $\| \cdot \|$ is the appropriate (unitarily invariant) norm on
$M(n)$ if $M=G$  (otherwise
$\| \cdot \|$ can be interpreted as a quotient norm on the corresponding
quotient of the
relevant Lie algebra;  cf. (14) in section 4 and comments following it). We
also point out that
all the  metrics we consider being bi-invariant,  any curve in $M$ can be
lifted
(by compactness and elementary properties of the $L_1$-norm)
to a ``transversal" curve in $G$ of the same length.
The correct abstract framework for these considerations is that of
Finsler geometry  (see e.g. \cite{busemann},
but the manifolds we consider being canonically embedded in normed spaces
we can afford to be more
``concrete".

For future reference we point out that the (operator) norm distance
and the corresponding intrinsic distance are related as follows
\begin{equation}
\|u-v\|=|1 - e^{i\,\rho(u,v)}|.
\end{equation}
\noindent This follows from Proposition 6 in the next section;
however, here we just wish to point out that the two metrics
differ by a factor of $\pi/2$ at the most and this particular fact is
implied by
the more or less obvious inequalities
$\rho(u,v) \geq \|u-v\| \geq |1 - e^{i\,\rho(u,v)}|$.
Since,  by definition, the corresponding
quotient metrics are distances between cosets,  the corresponding two metric
structure on homogeneous spaces are related in the way analogous to (6).
Accordingly, estimates of type (2) will
transfer easily from one metric to the other,  and the choice of the one
to work with will only be a matter of  convenience and/or elegance.

Because of the invariance of the metric under the action of $G$  ($=U(n)$
or $SO(n)$),
one can give estimates for the covering numbers of $M$ (analogous to those
of Lemma 1)
in terms of the Haar measure of balls (cf. the comment following (2)).
However,  since
the dependence of the measure of a ball on the radius is much less
transparent now than
in the ``linear" case,  such estimates are not necessarily useful.  To
overcome this
difficulty we ``linearize" the problem via the exponential map (composed
with the quotient map
$q:G \rightarrow M$ if necessary) and then use Lemma 2.  Since we operate
in the ``classical"
context,  the exponential map is the standard one
$$\exp x = e^x = \sum_{k =1}^{\infty} \frac{x^k}{k!} \mbox{ for }  x \in M(n),$$
and it will be normally sufficient to consider the restriction of $\exp$ to
$\cal{G}$ ($= u(n)$ or $so(n)$), the Lie algebra of $G$ ($=U(n)$ or $SO(n)$).
In order to be able to apply Lemma 2 we must
``understand"  $\Phi = q \, \circ \, \exp$;  specifically we need to know
for which
$K \subset \cal{G}$
we have  $\Phi(K) = M$ (or at least when $\Phi(K)$ is ``large") and for
which $K$
the restriction $\Phi_{|K}$ (resp. $\Phi_{| \Phi(K)}^{-1}$)  is Lipschitz.
Concerning the first point,
it is well known that, in our context, $\exp (\cal{G })$ $= G$.  Moreover,
we have

\medskip \noindent {\bf Lemma 3.} {\em Let }
$K = \{x \in {\cal{G}} \,: \, \|x\|_{\infty} \leq \pi \}$
{\em be the ball of radius } $\pi$ {\em in } $\cal{G}$ {\em in the operator
norm.  Then }
\newline {\em (a) } $\exp (K) = G$
\newline {\em (b) } $\exp$ {\em is one-to-one on the interior of } $K$.

\medskip The above is a special case of a more general fact for Lie groups,
 but in the present setting
can be seen directly from the fact that every unitary matrix can be
diagonalized,  with
the argument for $SO(n)$ being just slightly more complicated.

Lemma 3 asserts that $G$ resembles,  in a sense, a ball in $\cal{G}$.
However,  for our purposes
we need more quantitative information about $\exp$, which we collect in the
next lemma.

\medskip\noindent {\bf Lemma 4.} {\em For any unitarily invariant norm and
the corresponding
metric on } $G$ {\em  (extrinsic on intrinsic),  the map } $\exp : \cal{G}
$ $\rightarrow G$
{\em is a contraction. }
\newline {\em On the other hand, let } $\|\cdot\|$
{\em be any unitarily invariant norm and set, } for $\theta > 0$,
$$\phi(\theta) = \inf \{\frac{\|e^x - e^y\|}{\|x-y\|} \,: \, x, y \in
{\cal{G}}, \, x \neq y, \,
\|x\|_{\infty} \leq \theta, \, \|y\|_{\infty} \leq \theta \}.$$
{\em Then } $\phi(\theta) > 0$ {\em if } $\theta < \pi$. {\em  Moreover,  if }
$\theta \in [0,2\pi/3)$, {\em  then}
$$\phi(\theta) \geq \prod_{k=1}^{\infty} {(1 - |1 - e^{i \theta / 2^k}|)}.$$
{\em  In particular } $\phi(\theta) \geq .4$  {\em if } $\theta \leq \pi/4$.

\medskip\noindent {\bf Proof.} The first assertion is classical for the
extrinsic (norm)
metric and hence follows formally
for the intrinsic metric.  For the other assertions,  we observe first that
since
the derivative of the exponential map at $0$ is the identity,
\begin{equation}
\lim_{\theta \rightarrow 0^+} \phi(\theta) = 1.
\end{equation}
Let $x, y$ be as in the definition of $\phi(\theta)$.  We have
\begin{eqnarray*}
e^x - e^y &= &e^{\frac{x}{2}}(e^{\frac{x}{2}}-e^{\frac{y}{2}}) +
(e^{\frac{x}{2}}-e^{\frac{y}{2}}) e^{\frac{y}{2}} \\
&=&2(e^{\frac{x}{2}}-e^{\frac{y}{2}}) +
(e^{\frac{x}{2}}-I)(e^{\frac{x}{2}}-e^{\frac{y}{2}})
+(e^{\frac{x}{2}}-e^{\frac{y}{2}}) (e^{\frac{y}{2}}-I)
\end{eqnarray*}
and so,  by the ideal property of unitarily invariant norms,
\begin{eqnarray*}
\|e^x - e^y\| &= &\|e^{\frac{x}{2}}-e^{\frac{y}{2}}\| (2 -
\|e^{\frac{x}{2}}-I\| - \|e^{\frac{y}{2}}-
I\|) \\
&\geq&\phi(\frac{\theta}{2}) \|x-y\| \cdot (1 - |1 - e^{\frac{i \theta}{2}}|).
\end{eqnarray*}
\noindent Iterating and using (7) we obtain the third (and hence the last)
assertion of the lemma.  For the second assertion
($\phi(\theta) > 0$ if $\theta < \pi$,  not used in the sequel),   we just
briefly sketch
the argument for $G = U(n)$.
Let $\theta \in (0,\pi)$ and $\delta > 0$.  Consider first the case of the
operator norm.
 We need to show that if $A, B$ are
Hermitian with  spectra contained in $[-\theta, \theta]$ and $\|e^{iA} -
e^{iB}\| \leq \delta$,
then $\| A - B \| \leq C(\theta) \delta$,  where $C(\theta)$ depends only
on $\theta$
(and not on $A, B$ or $n$). By  \cite{bhatia}, Theorem 13.6, the
eigevalues of $A$ and $B$ (multiplicities counted) are,  ia a certain
precise sense, ``close",
and so,  by perturbation, we may assume that they are identical;  we may
also assume that
all those eigenvalues are inegral multiples of $\delta$. Let $u \in U(n)$
be such that $B = uAu^{-1}$;
we need to show that $\|e^{iA}u - ue^{iA}\| \leq 4\delta$ implies
$\|Au - uA\| \leq C(\theta) \delta$ and this follows by writing $u$ as a
``block matrix" in the
spectral subspaces of $A$.  For a general unitarily invariant norm we note
that the assertion is
roughly equivalent to uniform boundedness (with respect to the norm in
question and with a bound
depending only on $\theta$)
of the inverse of the  derivative of the exponential map.  That derivative
is, in a proper
orthonormal basis,  an antisymmetric ``Schur multiplier" (see \cite{va},
Theorem 2.14.3
and its proof).
As a consequence, the inverse is also a ``Schur multiplier" and so its norm
with respect to the
operator norm equals to the norm on the trace class $C_1$ (by duality) and
dominates the
norm with respect to any unitarily invariant norm by interpolation
(cf. \cite{tomczak}, \S 28 or \cite{gohberg}).

\medskip\noindent {\bf Remark.}  In all likelihood,  a version of Lemma 4
(and of Lemma 5
that follows) should be known,  at least for the operator norm,  but we
couldn't find a reference.
It would be nice to have an elegant proof which gives good constants in the
full range of
$\theta$ ($\in (0, \pi)$).  We point out that the
``Schur multiplier" argument indicated above provides a simple ``functional
calculus" proof
(with good constants) in the case of the Hilbert-Schmidt $C_2$-norm.

\medskip\noindent {\bf Lemma 5.} {\em Let } $G = U(n)$ {\em (resp.}
$SO(n)$) {\em and }
$\rho$ {\em the intrinsic metric on } $G$ {\em induced by a unitarily
invariant norm }
$\| \cdot \|$ {\em on } $M(n)$ {\em . Then,  for any } $x, y \in \cal{G}$,
$$\rho(e^{x+y}, e^x e^y) \leq \|[x,y]\|.$$

\medskip\noindent {\bf Proof.}  The argument is similar to,  but slightly
more complicated than that of
the previous lemma.  Denote,  for $t \geq 0$,
$$\psi(t) = max \{\rho(e^{t(x+y)}, e^{tx} e^{ty}),\rho(e^{t(x+y)}, e^{ty}
e^{tx})\}.$$
\noindent Clearly $\psi(0)=0$.  Moreover,  expanding the exponentials and
noting that
$\rho(u,v) / \|u-v\| \rightarrow 1$ as $\|u-v\| \rightarrow 0$  (this follows
from (6),  but can also be seen from the inequalities in the paragraph
following (6),
which do not depend on Proposition 6) we conclude that

\begin{equation}
\lim_{t \rightarrow 0} \frac{\psi(t)}{t^2} = \|[x,y]\|;
\end{equation}

\noindent Now
\begin{eqnarray*}
\rho(e^{x+y}, e^x e^y) & \leq \rho(e^{x+y}, e^{\frac{x+y}{2}}
e^{\frac{x}{2}} e^{\frac{y}{2}}) +
\rho(e^{\frac{x+y}{2}} e^{\frac{x}{2}} e^{\frac{y}{2}},
e^{\frac{x}{2}} e^{\frac{y}{2}} e^{\frac{x}{2}} e^{\frac{y}{2}}) \\
&+
\rho(e^{\frac{x}{2}} e^{\frac{y}{2}} e^{\frac{x}{2}} e^{\frac{y}{2}},
e^{\frac{x}{2}} e^{\frac{x+y}{2}} e^{\frac{y}{2}})  +
\rho(e^{\frac{x}{2}} e^{\frac{x+y}{2}} e^{\frac{y}{2}},e^x e^y) \\
& = 3\rho(e^{\frac{x+y}{2}}, e^{\frac{x}{2}} e^{\frac{y}{2}})
+ \rho(e^{\frac{y}{2}} e^{\frac{x}{2}}, e^{\frac{x+y}{2}})
\end{eqnarray*}

\noindent Hence $\psi(1) \leq 4 \psi(\frac{1}{2})$ and,  by the same argument,
$\psi(t) \leq 4 \psi(\frac{t}{2})$ or $\frac{\psi(t)}{t^2}
\leq \frac{\psi(\frac{t}{2})}{(\frac{t}{2})^2}$
for $t \geq 0$.  In combination with (8) this implies the lemma.

\section{SOME NON-RIEMANNIAN GEOMETRY}

Our last auxiliary result involves non-Riemannian geometry of $G$ = $U(n)$
(or $SO(n)$).
It is very well known (in a much more general context) that if $G$ is
endowed with
a bi-invariant Riemannian structure (which is, in our case, the one induced
by the
Hilbert-Schmidt $C_2$-norm on $M(n)$), then the geodesics of $G$ are
exactly the cosets of
one-parameter subgroups (see \cite{helgason}, p. 148, Ex. 5, 6). It is not
immediately
clear how general is this phenomenon.
Since geodesics are normally defined via affine connections, we should make
clear
here that we emphasize the ``metric" approach: a curve in a manifold $M$
endowed with a metric
is a geodesic if it {\em locally} realizes the (intrinsic) distance between
points as
explained in the paragraph containing (4)  (the argument can be presumably
rewritten by
starting from the affine connection induced by the group structure,
though). We have

\medskip\noindent {\bf Proposition 6.} {\em Let } $\| \cdot \|$
{\em be a unitarily invariant norm on }
$M(n)$ {\em and } $\rho$ {\em the induced intrinsic metric on }
$G$ ($ = U(n)$ {\em or } $SO(n)$). {\em Then {\em
\newline (a) {\em cosets of one parametric semigroups (i.e. curves of the form }
$\gamma (t) = ue^{tx}$, $u \in G$, $x \in \cal{G}$) {\em are geodesics in }
$(G, \rho)$
\newline (b) {\em if } $\| \cdot \|$ {\em is strictly convex (which happens
in particular if }
$\| \cdot \| = \| \cdot \|_p$ {\em for some } $p \in (1, \infty)${\em ),
then all geodesics are, up to a change of parameter, of the form given in
(a) (or arcs of curves of such form) }
\newline (c) {\em if, furthermore, the spectrum of } $u^{-1}v$
{\em does not contain } $-1${\em , the curve
of shortest length (geodesic arc) connecting } $u$ {\em and } $v$ {\em is
unique.}

\medskip {\bf Remarks.} (1) A unitarily invariant norm on $M(n)$ is
strictly convex
(i.e. the unit sphere $\{x \, : \, \| x \| = 1\}$
does not contain a segment) iff the associated
symmetric norm on ${\bf R}^n$ is (cf. \cite{gohberg}). The operator norm
and the
trace class $C_1$-norm are {\em not } strictly convex.
\newline (2) It is a somewhat delicate issue how smooth should be the
curves we consider,
particularly since the results from \cite{bhatia}, to which we refer quite
heavily, do not
fit {\em precisely } our needs {\em exactly} in that respect.
For the purpose of following the proof below,  the reader should
think of all functions as being at least $C^1$. As indicated in the
previous section,  absolutely
continuous functions provide a convenient framework;  we comment on the
fine points of
the present argument and on their relevance to \cite{bhatia} at the end of
the proof.

\medskip \noindent{\bf Proof.} Since $SO(n)$ is a Lie subgroup of $U(n)$,
it is enough to prove the lemma for the latter.
\noindent (a) By unitary invariance of the metric it is enough to prove the
assertion for
``short arcs of one parameter semigroups",  i.e. curves of the form
$\gamma_0 (t) = e^{tx}$, $t \in [0, 1]$ where $\| x \|$ is ``small".
By  (5),  the length $\ell (\gamma_0 )$ equals $\| x \|$.  We need to show
that any curve $\gamma$
in $U(n)$ connecting $\gamma_0 (0) = I$ and $\gamma_0 (1) = e^{x}$ is of
length $\geq \| x \|$.
The argument is based on  (and in fact very close to) the results on spectral
variation of unitary matrices presented in \cite{bhatia}, \S 13, 14.
Indeed,  Theorem 14.3 and
Remark 14.4 of \cite{bhatia}, when specified to unitary matrices, say in
effect that the map
$\Sigma$ associating
to a matrix its spectrum is a contraction when considered as acting from
$(U(n), \rho)$ to
$({\bf C}^n, \| \cdot \|)/S_n$,  where $\| \cdot \|$ is now the symmetric norm
associated to the unitarily invariant norm in question (a priori defined on
${\bf R}^n$, but
canonically extendable to ${\bf C}^n$) and $S_n$ is the symmetric group
acting on ${\bf C}^n$ by
permuting the coordinates (we recall that we are working with the {\em
intrinsic} metric).
This in turn implies that, for any curve $\gamma$ in $U(n)$ connecting $I$
and $e^{x}$, we have
\begin{equation}
\ell (\gamma) \geq \ell (\Sigma(\gamma)).
\end{equation}
  Since $\gamma$ lies in $U(n)$,
its image  under $\Sigma$ lies in ${\bf T}^n/S_n$,  where
${\bf T} = \{z \in {\bf C} \,: \, |z| = 1 \}$.  Let now
$E \,: \, {\bf R}^n \rightarrow {\bf T}^n$ be the exponential map of the
group ${\bf T}^n$
(i.e., $E((\xi_k)_{k=1}^n) = (\exp(i \xi_k))_{k=1}^n$) considered as acting
from
$({\bf R}^n, \| \cdot \|)$ to ${\bf T}^n$ equipped with the intrinsic
metric inherited from
$({\bf C}^n, \| \cdot \|)$.  Then $E$ is a local isometry,  e.g. it is an
isometry when
restricted to $[-\pi/2, \pi/2]^n$.  Consequently,  the appropriately
restricted induced map
$\tilde{E} \,: \, {\bf R}^n / S_n \rightarrow {\bf T}^n/S_n$ is also an
isometry and so,
if  $\| x \|$ is ``small",
the length of $\gamma_1 = \tilde{E}^{-1} (\Sigma(\gamma))$ is the same as
that of
$\Sigma(\gamma)$,  in particular,  by (9), $\ell (\gamma_1) \leq \ell (\gamma)$.

Now observe that $\gamma_1$ connects $0$ ($\in {\bf R}^n$)  and
$\lambda = (\lambda_1, \dots, \lambda_n)$,  where $(i \lambda_1, \dots, i
\lambda_n)$ are
eigenvalues of $x$ with multiplicities.  Clearly  $\| x \| = \| i \lambda
\| = \| \lambda \|$
and they are all equal to the distance between $0$ and $\lambda$ in ${\bf
R}^n/S_n$.  Accordingly
\begin{equation}
\| x \| \leq  \ell (\gamma_1) \mbox{ (in } {\bf R}^n/S_n \mbox{) }
= \ell (\Sigma(\gamma)) \mbox{ (in } {\bf T}^n/S_n \mbox{) }
\leq  \ell (\gamma) \mbox{ (in } (U(n),\rho) \mbox{)}
\end{equation}
and the assertion (a) is proved.

(b) Again,  it is enough to consider curves of the form
$\gamma_0 (t) = e^{tx}$, $t \in [0, 1]$ where $\| x \|$ is ``small".  We
have to
show that if $\gamma$ is a curve in $U(n)$ connecting $I$ and $e^{x}$ such
$\ell (\gamma) = \| x \| = \ell (\gamma_0)$,  then (after a change of
parameter, if
necessary), $\gamma = \gamma_0$. In other words,  we need to investigate
the cases of
equality in (10).
\newline Concerning the first inequality in (10), we observe that since
$\| \cdot \|$ is strictly convex, a
straight line segment in $({\bf R}^n,\| \cdot \|)$ with endpoints $\lambda,
\lambda'$
is {\em strictly} shorter than any other curve connecting these points.
As a consequence, the same is true
for the corresponding quotient metric in ${\bf R}^n/S_n$ provided the
distance between
$\lambda, \lambda'$ in the quotient metric equals $\| \lambda - \lambda'
\|$ (in particular
if one of the endpoints is $0$).
Thus we may have equality in the first inequality in (10) only if
$\gamma_1$ is a segment,
i.e. if (after a change of parameter, if
necessary) the spectrum of $\gamma (t)$, or $\Sigma (\gamma) (t)$, equals
$(\exp (i \lambda_1 t), \dots, \exp (i \lambda_n t))$,  where,  as before,
$(i \lambda_1, \dots, i \lambda_n)$ are eigenvalues of $x$.  Let $i  \mu_1,
\dots, i \mu_m$ be all
{\em distinct} eigenvalues of $x$,  and let,
for $k=1,\dots,m$ and $t \in [0,1]$,
$P_k (t)$ be the spectral projection of $\gamma (t)$ corresponding to the
eigenvalue
$\nu_k (t) = \exp (i \mu_k t)$  (since $x$ was assumed to be ``small",
$\nu_k$'s are also distinct).
If follows from elementary functional calculus that $P_k (t)$ is a
continuous function of $t$ and,  moreover, has the same smoothness (in $t$)
as $\gamma (t)$.
To prove that $\gamma = \gamma_0$ we need to show that  $P_k (t)$ is
constant in $t$ for
$k=1,\dots,m$.  To this end, we need to analyze the equality case in the
second inequality
in (10), and for that we must go into the proof of Theorem 14.3 of
\cite{bhatia}.  In our setting
and notation, it is shown there (p. 68, equation (14.5)) that
\begin{equation}
\ell (\Sigma(\gamma)) \leq \int_0^1 \|P_{\gamma(t)}\gamma'(t) \| dt
\leq \int_0^1 \|\gamma'(t) \| dt = \ell (\gamma),
\end{equation}
where,  for $u \in U(n)$,  $P_u$ is a (necessarilly contractive) orthogonal
projection in
$M(n)$ onto the commutant of $u$ (in our case only the restriction of the
projection to $\cal{G}$
is relevant).  By the strict convexity of $\| \cdot \|$,  we have
$\|P_u x \| < \| x \|$ unless $x$ is in the
range of $P_u$  and so the second inequality in (11)
(and hence the second inequality in (10)) is strict unless $[\gamma'(t),\gamma(t)] = 0$
for amost all $t \in [0,1]$.  Now, $\gamma(t)$ being normal, $\gamma'(t)$
must also commute
with the spectral projections of $\gamma(t)$,  i.e.
\begin{equation}
[\gamma'(t),P_k (t)] = 0  \mbox{ for } k=1, \dots ,m.
\end{equation}
Since $\gamma(t) = \sum_{j=1}^m \nu_j (t) P_j (t)$,  (11) translates into
\begin{equation}
[\sum_{j=1}^m \nu_j (t) P_j' (t), P_k (t)] = 0  \mbox{ for } k=1, \dots ,m.
\end{equation}
For (small) $h>0$, choose $v(h) \in U(n)$ so that it conjugates
$\sum_{j=1}^m \nu_j (t) P_j (t)$  and
$\sum_{j=1}^m \nu_j (t) P_j (t+h)$  (it then also conjugates $P_k (t)$ and
$P_k (t+h)$ for
$k=1, \dots ,m$)
and so that $v(0) = I$ and $v(h)$ depends ``smoothly" on $h$.
Such a ``smooth" choice is possible wherever $\gamma$ is ``smooth":  as we
indicated earlier,
$(P_k)_{k=1}^m$ is
then also a ``smooth" curve in the the "generalized Stieffel manifold"
(i.e., the quotient
of $U(n)$ by the commutant of $\gamma(t)$; see section 5 for an
elaboration of the concept and its framework) and so it can be locally
lifted to an
``equally smooth" curve in $U(n)$.  If $y = y(t) = v'(0)$,  then it is
easily seen that
\begin{equation}
P_k' (t) = [y, P_k (t)]  \mbox{ for } k=1, \dots ,m
\end{equation}
and so (13) translates into
$$[y, \sum_{j=1}^m \nu_j (t) P_j (t)] \in \mbox{ commutant of } \gamma (t),$$
and so (dropping $t$,  which is fixed, from the formulae),
$$[y, \sum_{j=1}^m \nu_j P_j] = \sum_{k=1}^m P_k [y, \sum_{j=1}^m \nu_j
P_j] P_k ,$$
which is easily seen to be $0$.  On the other hand,
$$[y, \sum_{j=1}^m \nu_j P_j] = \sum_{k=1}^m \sum_{j=1}^m (\nu_j - \nu_k)
P_k y P_j.$$
Since $\{ P_k M(n) P_j,\;  k=1, \dots ,m,\, j=1, \dots ,m \}$
form an orthogonal decomposition of $M(n)$,
the above double sum can be $0$ only if all the terms are $0$.  Recalling
that $\nu_k$ are
distinct,  we deduce that  $j \neq k$ implies $P_k y P_j = 0$.
Hence $[y, P_k (t)] = \sum_{j \neq k}  P_j y P_k - P_k y P_j = 0$ and so,
by (14), $P_k'(t) = 0$
for $k=1, \dots ,m$.
Since $t$ was an arbitrary point in $[0,1]$ at which $\gamma (t)$ was
``smooth",  it follows that
$P_k'(t) = 0$ for almost all $t$ and so each $P_k (t)$ is in fact constant,
as required.
This proves part (b).
\newline (c) This is immediate;  by (b), the geodesic must be (up to a
change of parameter)
of the form $\gamma (t) = ue^{tx},  t \in [0,1]$ with $\ell (\gamma ) = \|
x \|$ and
$e^x = u^{-1}v$.  Since the spectrum of $u^{-1}v$ does not contain $-1$,
the last equation
has unique solution $x_0$ such that $\| x_0 \|_{\infty} < \pi$,  and any
other solution $x_1$
verifies $\| x_1 \| > \| x_0 \|$,  thus leading to a {\em strictly } longer
curve.

\medskip \noindent{ {\bf Comments.} (1) We need to clarify that the results
presented in \cite{bhatia}
and used above were obtained under ``piecewise $C^1$" assumptions,  both
for the curve and
for (roughly speaking and in our notation) the distance function in
$({\bf C}^n, \| \cdot \|)/S_n$ applied to $\Sigma(\gamma)$.  Such framework
was good enough for
\cite{bhatia}  (their analysis, even though it fails to identify {\em all}
singular points
of the latter ``distance function" in the general case, can be easily
patched by approximating a
general unitarily invariant norm by a smooth one),  but is insufficient in
our context,
particularly for settling the equality cases. As indicated before,  one
solution is to consider
absolutely continuous curves to insure that variation of a function can be
expressed in terms
of its derivative;  we tacitly used that property several times in our
argument).
Another ``fine point" that makes the argument work is the fact that the
``distance function in
$({\bf C}^n, \| \cdot \|)/S_n$",  being an minimum of a finite number of
convex functions,
is differentiable almost everywhere (also when composed with an absolutely
continuous curve)
and, moreover,  it  has directional derivatives {\em everywhere}.
\newline (2) In the last part of the argument we did show that,  under certain
additional assumptions,  if a curve $\gamma$ in $U(n)$ verifies
$[\gamma(t), \gamma '(t)] = 0$,  then it must be contained in a commutative
subgroup.
It appears  likely that the argument carries over to a more general class
of ``sufficiently
smooth" curves in $U(n)$:  what we seem to be exploiting is that the patern of
equatities between eigenvalues is constant.

\section{THE ``1-COMPLEMENTED" SUBGROUPS}

Let $G = SO(m)$ or $U(n)$ and $\cal{G}$ ($= so(n)$ or $u(n)$) - the Lie
algebra of $G$.
Let $H$ be a connected Lie subgroup of $G$, $\cal{H}$ - the corresponding
Lie subalgebra of $\cal{G}$ and  $M = G/H$.  In this
section we concentrate on the case when $G$ and $M$ are endowed with metric
structures induced by the operator norm,  we will call the respective
metrics by $\rho$ and $\rho_M$  (for technical purposes, we may use other
unitary
ideal norms, though).  The purpose of this section is to prove,  in the
above context,
an estimate of type (2) for an abstract class of homogeneous spaces that
contains
$U(n)$, $SO(n)$ and the Grassmannians $G_{n,k}$.
The argument will depend on a careful analysis of the
exponential map $\exp \, : \, {\cal{G}} \rightarrow G$ and maps obtained
from it;  as an
illustration we point out here that the following result from \cite{Iowa} is an
immediate consequence of the results from the preceding two sections.

\medskip\noindent {\bf Theorem 7.} {\em If } $G = SO(n)$ {\em or } $U(n)$
{\em (endowed with the operator norm or the induced intrinsic metric }
$\rho${\em)  and } $\varepsilon \in (0, 2]${\em , then }
$$(\frac {c}{\varepsilon} )^d \leq N(G, \varepsilon) \leq (\frac
{C}{\varepsilon} )^d,$$
\noindent {\em where } $d$ {\em is the (real) dimension of }$G$  {\em and }
$c$ {\em and } $C$ {\em are universal numerical constants.}

\medskip\noindent {\bf Proof.} By Lemma 3 and the first assertion of Lemma 4,
$\exp$ is a contractive surjective map
from the (closed) ball of radius $\pi$ in $\cal{G}$ to $G$.  Consequently,
Lemma 2 applied with
$\Phi = \ exp$ and $L = 1$ and combined with the second inequality in (1)
(or,  more precisely,
the upper stimate on $N''$ given by Lemma 1)
yields the upper estimate for  $N(G, \varepsilon)$.  The lower estimate is
obtained
similarly by applying Lemma 2 to $\Phi = \ exp^{-1}$,  $L = 2.5$ and
$K = \{u \in G \, : \, \rho(u, I) \leq \pi/4\}$, using the first inequality
in (1)
(or, again more precisely,
the lower estimate on $\tilde{N}$ given by Lemma 1),
the last assertion of Lemma 4 and Proposition 6(a).
(Proposition 6(a) is needed only for the definition of $K$ and its use may
be avoided here.)

\medskip If $M \neq G$,  the approach will be similar,  but the situation
is (necessarily) more
complicated.  Let $q \,: \, G \rightarrow G/H = M$ be the quotient map and
consider
the short exact sequence
$0 \rightarrow H \rightarrow G \rightarrow M \rightarrow 0$ and the induced
sequence of
maps between the tangent spaces (at resp. $I \in H$, $I \in G$ and $H \in
M=G/H$)
\begin{equation}
0 \rightarrow {\cal{H}} \rightarrow {\cal{G}} \rightarrow T_H M \rightarrow 0
\end{equation}
and so $T_H M$ can be identified with the quotient space $\cal{G}/\cal{H}$;
we mean by that
{\em isometrically} identified whenever all metric structures are induced
by a given
unitarily invariant norm. Since the derivative of the exponential map at
$0$ is the identity
(in particular an isometry),  we can realize that identification by the
canonical
factorization of the derivative of $q \, \circ \, \exp $ at $0$ (which maps
$\cal{G}$ to
$T_H M$ and vanishes on $\cal{H}$) through $\cal{G}/\cal{H}$.  This shows
that (at least small)
neighbourhoods in $M$ resemble balls in the normed space $\cal{G}/\cal{H}$
and gives
some heuristic evidence that inequalities of type (2) may hold for $M$.
However,  for a proof of such an inequality one needs
``uniform isomorphic" (rather than ``infinitesimal") estimates,  and we
will obtain these under
some additional technical assumptions.  Since the additive structure on
$\cal{G}$ and the group
structure on $G$ are not intertwined by the exponential (or any other) map,  it will be
more convenient to identify $\cal{G}/\cal{H}$ with $\cal{X} = \cal{H}^{\perp}$
(the orthogonal complement of $\cal{H}$ in $\cal{G}$) and to consider
$\Phi = q \, \circ \, \exp_{|\cal{X}}$,  hoping that the direct sum
$\cal{X} \oplus \cal{H} = \cal{G}$ is ``well-behaving" with respect to the
operator norm
(or any other unitarily invariant norm that we may need to consider),  which
happens in many natural examples.
This leads to our first invariant related to a homogeneous space.  We set

\begin{equation}
\kappa(M) = \| P_{\cal{X}} \| = \| I - P_{\cal{H}} \|,
\end{equation}

\noindent where $P_E$ denotes the orthogonal projection onto $E$ and $\|
\cdot \|$ is calculated
with respect to the operator norm.

Before stating the results, we will present a simple but illuminating
example which
shows that,  in general,  the ``linearization" of $M$ of the type suggested
above may work
{\em only} on the ``infinitesimal" scale (i.e. only very small neighbourhoods
are ``equivalent" to balls in the tangent space),
and which leads to one more invariant of $M$.  Let  $G = U(n)$ and $H =
SU(n)$.
It is then easily seen that $M = U(n)/SU(n)$ is {\em isometric } to a
circle of radius $1/n$ and
so covering numbers of $N(M,\varepsilon)$ are ``trivial" if $\varepsilon >
\pi/n$.
(Since $M$ is 1-dimensional,  it is necessarilly ``isotropic" and so there are
neighbourhoods resembling segments of size comparable to the diameter
of $M$; in particular (2) still holds. However,  one can also produce
``nonisotropic"
examples: consider e.g.,  $H = \{ I \} \times SU(n-1) \subset U(n) = G$.)
The reason for this phenomenon
is that $SU(n)$ (or ($\cal{H}$,  via the exponential map) is very ``densely
woven"
into $U(n)$.  For example,  $e^{2 \pi i/n}I \in SU(n)$ and
$\| e^{2 \pi i/n}I - I \| < 2 \pi /n$  (more precisely,
$\rho (e^{2 \pi i/n}I, I) = 2 \pi /n$ by
Proposition 6(a)),  even though the shortest path connecting $I$ and $e^{2
\pi i/n}I$ and
contained in $SU(n)$ is of length $2 \pi (1 -1/n)$ (this follows from the
proof of
Proposition 6(a), the length in question must be
$\geq$ than the length of the shortest path connecting
$(-2 \pi + 2 \pi/n, 2 \pi/n, \dots,2 \pi/n)$ and $0$ in $\ell_{\infty}^n / S_n$
that is contained in the plane $\{ (x_k) \in {\bf R}^n \, : \, \sum x_k = 0
\}$;
another way to express this is that $e^{2 \pi i/n}I = e^x$ with $x \in \cal{H}$
forces $\| x \| \geq 2 \pi (1 -1/n)$).
To quantify the phenomenon we introduce the following concept.  Given
$\theta > 0$,
we will say that  a closed connected Lie subgroup $H$ of $G = U(n)$
(or $SO(n)$) is $\theta$-woven if whenever
$u \in H$ satisfies $\rho (u,I) \leq \theta$  ($\rho$ is the intrinsic
metric induced by the operator norm $\| \cdot \|_{\infty}$),  then there exists
$x \in \cal{H}$, $\| x \|_{\infty} < \pi$ such that $u = e^x$. If $M =
G/H$, we set

\begin{equation}
\theta(M) = \sup \{ \theta > 0 \, : \, H \mbox{ is } \theta \mbox{-woven} \} =
 dist(I,H \backslash \exp (B_{\cal{H}}(\pi)),
\end{equation}

\noindent the distance being calculated using $\rho$.  We then have

\medskip\noindent {\bf Theorem 8.} {\em In the notation above,  assume that }
$\kappa(M) = 1$. {\em  Then, for any } $\varepsilon \in (0,diam\, M]$,
$$N(M, \varepsilon) \leq (\frac {C diam\, M}{\varepsilon })^d,$$
\noindent {\em where } $d$ {\em is the (real) dimension of } $M$,  $diam\,
M$ {\em is
calculated with respect to } $\rho_M$,  {\em and } $C > 0$ {\em is a
universal constant.
Moreover,  if } $\varepsilon \in (0,\theta(M)/4]$, {\em  then }
$$N(M, \varepsilon) \geq (\frac {c \theta (M)}{\varepsilon })^d ,$$
\noindent {\em where } $c > 0$ {\em is a universal constant. The last estimate
holds also if } $\kappa(M) > 1$, {\em  but the constant $c$ may then depend
on } $\kappa(M)$.

\medskip\noindent {\bf Proof.} As suggested earlier,  the proof will
involve applying
Lemma 2 to the (properly restricted) map $q \, \circ \, \exp_{|\cal{X}}$
and its inverse,
where $\cal{X}$ is the orthogonal complement in $\cal{G}$ of $\cal{H}$  (the Lie
subalgebra of $\cal{G}$ corresponding to the subgroup $H$),  and will be
based on two
lemmas that follow. Given $r > 0$,  let $B_{\cal{X}}(r)$ be the ball in
$\cal{X}$
of radius $r$ (with respect to the operator norm) and centered at the origin.
We then have,  in the notation of Theorem 8:

\medskip\noindent {\bf Lemma 9.}  {\em If } $\kappa(M) = 1$, {\em  then }
$q(exp(B_{\cal{X}}(diam\, M))) = M$.

\medskip\noindent {\bf Lemma 10.} {\em There exist positive constants }
$\lambda = \lambda (\kappa(M))$ {\em and }
$r_0 = r_0(\kappa(M))$ {\em such that if } $r = \min \{r_0,\theta (M)/4 \}$
{\em and }
$x, x' \in B_{\cal{X}}(r)$, {\em  then }
$$\rho_M (q(e^x), q(e^{x'})) \geq \lambda \| x - x' \|.$$

\medskip\noindent {\bf Remarks.} (i) Calculating $\theta (M)$ is not
difficult,
particularly when
$H$ is semisimple.  Indeed,  suppose $\theta (M) < \pi$  (clearly the maximal
possible value)
and let $u \in H \backslash \exp (B_{\cal{H}}(\pi))$ be such that $u = e^h$,
$\| h \|_{\infty} \geq \pi$, while $\| u - I \|_{\infty} = \theta(M) < \pi$,
in particular $u = e^x$ for some $x \in \cal{G}$, $\| x \|_{\infty} < \pi$.
Since the commutants of $u$ and $x$ are the same, it follows that $\theta
(M)$ is
``witnessed" inside a torus $\cal{T}$ in $G$ containing $u$;  moreover,
$\cal{T}$
may be assumed to contain the one-parameter semigroup
$\{ e^{tx} \, : \, t \in {\bf R} \}$ and to be such that $\cal{T} \cap H$ is
maximal in $H$.  Consequently,  to determine $\theta (M)$ we only need to
examine
maximal tori in $H$ and their extensions to maximal tori in $G$.
This is particularly easy if $H$ is semi-simple:  all configurations of the
tori in question are then
related by conjugation,  and since the metric we consider is invariant under
conjugation,  it suffices to check just one such configuration. Such an
examination
will also reveal that $\theta(M) = \pi$,  should that be the case.
\newline (ii) Since $diam\, U(n) = \pi$ and $q$,  being a quotient map,
is a contraction,  one {\em always} has $diam\, M \leq \pi$. In any case,
by Proposition 6,
one can always calculate $diam\, M$ by examining images of one-parameter
semigroups
of $G$ under $q$.
\newline (iii) If $M = G_{n,k}$  (the Grassmann manifold), one verifies
directly that
$\kappa(M) = 1$ and  $diam\, M = \pi/2$.  The former follows from the fact
that $\cal{X}$ consists of those matrices in $\cal{G}$ ($ = u(n)$ or $so(n)$)
that are of the form  (cf. (3))
$$\left [
\begin{array}{cc}
0 & x  \\
-x^* & 0
\end{array}
\right ] $$
The latter is elementary
and presumably known:  for two ($k$-dimensional) subspaces $E, F$ of
${\bf R}^n$ (resp. ${\bf C}^n)$),  $\rho_{G_{n,k}}(E,F)$ is the largest of
the main angles between $E$ and $F$ (see (4) and comments following it
for the framework and
e.g., \cite{Iowa},  p. 174 for the more precise analysis). Finally,  it follows
immediately from the previous remark that $\theta(G_{n,k}) = \pi$ (the
maximal tori in $H$
are also maximal in $G$).
\newline (iv) If $\kappa(M) = 1$, one can take in Lemma 10 $r = .12$ and
$\lambda = .4$ and if,
moreover, $x'=0$, one can take $r = 5/9$ and $\lambda = .4$;  it follows that
$q \, \circ \, \exp \,(B_{\cal{X}}(5/9)) \supset \{U \in G/H \, : \,
\rho(U, q(I)) \leq 2/9 \}$
and that $q \, \circ \, \exp ^ {-1}$ restricted to any of these two sets is
Lipschitz
with constant $4.5$.
\newline (v) The proof gives $\lambda (t)$ and $r_0(t)$ to be of order  $1/t$.
The argument would be slightly more efficient if we considered $\cal{X}$ as
endowed with the quotient norm $\cal{G}/ \cal{H}$,
which is more natural in the context.

Assuming the two lemmas above,  Theorem 8 is shown almost exactly as
Theorem 7:
one applies Lemma 2,  first with $\Phi = q \, \circ \, \exp_{|\cal{X}}$, $L
= 1$ and
$K = B_{\cal{X}}(diam\, M)$ for the upper estimate and then with the inverse
map restricted to $K = q(exp(B_{\cal{X}}(r)))$ and with $L = \lambda^{-1}$
for the lower estimate. (All the fine points are hidden in Lemmas 9 and 10.)

\medskip\noindent {\bf Proof of Lemma 9.}  We will show that,  for every  $p \in [2, \infty)$,  setting
$$K_p = \{x \in {\cal{X}} \, :\, \|x\|_p \leq n^{1/p}  diam \, M \}$$
(i.e. $K_p$  is a ball in $\cal{X}$ of radius $n^{1/p} diam \, M$ in the
Schatten $C_p$-norm
$\| \cdot \|_p$), we have
\begin{equation}
q(\exp (K_p)) = M;
\end{equation}
the assertion of the Lemma will then follow by letting $p \rightarrow
\infty$.  To this end,
observe that since the {\bf R}-linear orthogonal projection $M(n)
\rightarrow \cal{G}$
is of norm one  (with respect to {\em any} unitarily invariant norm),
$\kappa (M)$ equals
to the norm of the orthogonal projection from $M(n)$ onto $\cal{X}$.  Now,
since
$(M(n), \| \cdot \|_p)$ is a complex interpolation space between
$(M(n), \| \cdot \|_{\infty})$  and $(M(n), \| \cdot \|_2)$,  it follows
that $P_{\cal{X}}$
is also contractive with respect to the $C_p$-norm  (more generally,  of norm
$\leq \kappa(M)^{1 - 2/p}$).  Furthermore,  since the  $C_p$-norm is
strictly convex
for $p \in (1, \infty)$,  we conclude that
\begin{equation}
y \notin {\cal{X}} \Rightarrow \| y - P_{\cal{H}} y \|_p = \| P_{\cal{X}} y
\|_p < \| y \|_p .
\end{equation}
For clarity,  we will denote by $M_p$ the manifold $M$ equipped with the
quotient metric
$\rho_{p,M}$ induced by
the Schatten $C_p$-norm. Note that since the operator norm and the $C_p$-norm
differ by factor $n^{1/p}$ at the most,  we have $diam\, M_p \leq
n^{1/p}diam\, M$.
Let $gH \in M_p$ and let $\gamma$ be the shortest geodesic in $M_p$ connecting
$H$ and $gH$,  then $\ell (\gamma) \leq n^{1/p}diam\, M$.  Let
$\tilde{\gamma}$ be a
transversal lifting of $\gamma$ to $G$,  i.e. a
curve in $G$ such that $q \, \circ \, \tilde{\gamma} = \gamma$ and
$\ell (\tilde{\gamma}) = \ell (\gamma)$.  Then of course $\tilde{\gamma}$
is a geodesic
in $G$ (with respect to the intrinsic metric $\rho_p$ induced by the
$C_p$-norm)
and without loss of generality we may assume that the initial point of
$\tilde{\gamma}$ is $I$.  By
Proposition 6,  $\tilde{\gamma}$ must be (perhaps after a change of
parameter) of the form
$\tilde{\gamma}(t) = e^{ty}, \, 0 \leq t \leq 1$ for some $y \in \cal{G}$ and
$\ell (\tilde{\gamma}) = \| y \|_p \leq n^{1/p}diam\, M$, and so (18) will
follow if we
show that $y \in \cal{X}$.  Indeed,  if that was not the case, (19) would
imply that
$\| y - P_{\cal{H}} y \|_p < | y \|_p$ and so,  for $t > 0$ sufficiently
small we
would have
$$\rho_{p,M}(e^{ty} H, H) \leq \rho_p(e^{ty},  e^{tP_{\cal{H}} y})
\leq \| ty - tP_{\cal{H}} y \|_p < t \| y \|_p$$
and consequently
$$ \begin{array}{ll}
\ell (\tilde{\gamma}) & = \rho_{p,M}(H,gH) = \rho_{p,M}(H,e^y H) \\
& \leq \rho_{p,M}(H,e^{ty} H) + \rho_{p,M}(e^{ty} H,e^y H)
< t \| y \|_p + (1-t) \| y \|_p \; ,
\end{array} $$
a contradiction.  This proves Lemma 9.

\medskip\noindent {\bf Proof of Lemma 10.}  We need to show that if
$\; x, x' \in B_{\cal{X}}(r)$ and $h \in \cal{H}$, then
$$\Delta \equiv \rho (e^{x'}, e^xe^h) \geq \lambda \| x - x' \|.$$
Since $\| x - x' \| \geq \rho (e^{x'}, e^x) = \rho (e^{-x} e^{x'},I)$
and $\Delta = \rho (e^{-x} e^{x'},e^h)$,  it is enough to consider $h \in
\cal{H}$
such that
$$\rho (e^h,I) \leq (1 + \lambda) \| x - x' \| \leq (1 + \lambda) 2r \leq 4r.$$
If $r \leq \theta (M)/4$  (or just  $2(1 + \lambda)r \leq \theta (M)$),  it
follows
from the definition of $\theta (M)$ (i.e. (17)) that $h \in \cal{H}$ may be
further assumed to satisfy $\| h \|_{\infty} < \pi$,  hence
$$\| h \|_{\infty} = \rho (e^h,I) \leq (1 + \lambda) \| x - x' \|
\leq (1 + \lambda) 2r \leq 4r.$$
Now,  by Lemma 4 and Lemma 5,
$$ \begin{array}{ll}
\Delta  &\equiv  \rho (e^{x'}, e^{x}e^h) \geq \rho (e^{x'}e^{-h/2},
e^{x}e^{h/2})  \\
&\geq  \rho (e^{x' - h/2}, e^{x + h/2}) - \| [x', h/2] \| - \| [x, h/2] \| \\
&\geq  \phi (r + \| h \|/2) \| x - x' -h \| - 2r \| h \|\\
&\geq \phi (r + \| h \|/2) \| x - x' \| - 2r \| h \| \\
&\geq  (\phi (r + (1 + \lambda)r) - 2r(1 + \lambda))\| x - x' \| \\
&\geq (\phi (3r) - 4r) \| x - x' \|,
\end{array} $$
\noindent where $\phi ( \cdot )$ is the function from Lemma 4.
It is now clear from Lemma 4 that if
$r > 0$ is small enough,  then $\phi (3r) - 4r > 0$.  A more careful
calculation along the same lines shows that if $r = .12$,  then  $\lambda =
.4$ works
(as indicated in Remark (iv) above).
\newline Finally,  if $\kappa (M) > 1$,  we can only use
$\| x - x' -h \| \geq \kappa (M)^{-1} \| x - x' \|$ in the fourth
inequality in the
preceding argument.  This results in the last expression being
$(\phi (3r) \kappa (M)^{-1} - 4r) \| x - x' \|$,  which yields the assertion
of the Lemma with $r_0$ and $\lambda$ being of order $\kappa (M)^{-1}$.

\section{EXTENSIONS AND OTHER TRICKS.}

The scheme presented in the preceding section yields resonable estimates for
covering numbers $N(M, \varepsilon)$ (with respect to the metric induced by
the operator norm) of a homogeneous space $M = G/M$  whenever
$\varepsilon \leq diam\, M$ or $\varepsilon \leq \theta(M)$  (for the upper
and lower estimate respectively) and whenever $\kappa(M)$ is
appropriately controlled.  This leaves several cases and gray areas that
are not covered.
\newline (i) The range $\theta(M) < \varepsilon < diam\, M$ even if
$\kappa(M) = 1$ and,
in general, a clarification of the role of the ratio $diam\, M / \theta(M)$
(the lower and upper estimates differing roughly by $(diam\, M / \theta(M))^d$).
\newline (ii) The upper estimate whenever $\kappa(M) > 1$,  but still
"under control".
\newline (iii) The case when we do not control $\kappa(M)$.

With regards to (i), a modification of the example that led to the
definition of $\theta(M)$
suggested there ($H = \{ I \} \times SU(n-1) \subset U(n) = G$) shows that
it is
possible for $diam\, M$ and $\theta(M)$ to differ by a large factor (of
order $n$
in that case).  Even though an analysis of such cases is imaginably possible,
it would be clearly combinatorial and/or algebraic in nature and we do not
attempt it here.

Concerning (iii),  it is also conceivable that the phenomenon of having
$\kappa(M)$
``large" can be ``dissected" and expressed in terms of combinatorial/algebraic
invariants suggested above,  but,  again,  in the examples motivating this work
(see below) we have $\| P_{\cal{H}} \| = 1$ and hence $\kappa(M) \leq 2$.

This leaves the gap related to (ii):  the examples with, say,  $1 <
\kappa(M) \leq 2$
do naturally occur and it would be nice to have, at least for that case, an
upper
estimate for covering numbers of $M$ of the type $(C diam \, M/
\varepsilon)^d$.
Unfortunately,  we do not know how to settle that question in full generality.
Instead,  we present a ``bag of tricks" that allow to handle various
special cases.
This,and some comments concering covering numbers relative to metrics generated
by unitarily invariant norms other than the operator norm constitutes this
section.

The first observation is that trying to mimmick the proof of Lemma 9 in the
case when
$\kappa(M) > 1$ one arrives at the following picture. Let
$Q : \cal{G} \rightarrow \cal{G}/ \cal{H}$ be the quotient map, and
consider the
semi-norm $p$ on $\cal{G}$ defined by $p(x) = \| Qx \|_{\infty}$.
Let $\Lambda : \cal{G}/ \cal{H} \rightarrow \cal{G}$ be
a norm-preserving lifting of $Q$ (in general nolinear). The argument
mimmicking the
proof of Lemma 9 connects then geodesics in $M$ with ``rays" in the range
of $\Lambda$
and we could give upper estimates for entropy of $M$ if we were able to control
entropy of the range of $M$ (e.g. with respect to the semi-norm $p$).

The two specific subgroups of $U(n)$, for which estimates for covering
numbers of the
respective homogeneous spaces are of interest from the point of view of free
probability  (cf. \cite{voiculescu}, Remark 7.2 and \cite{voiculescupc}),
consist of unitaries of some
$C^*$-subalgebras of $M(n)$,  namely

\noindent (1) The ``block-diagonal" algebra:  the commutant of $\{ P_1,
P_2, \dots, P_m \}$,
where $P_j$'s are orthogonal projections
whose ranges form an orthogonal decomposition of ${\bf C}^n$.
\newline (2) The ``tensor-factor" algebra:  if  $n = mk$,  identify ${\bf
C}^n$ with
${\bf C}^m \otimes{\bf C}^k$ and consider matrices of the form $I \otimes x$, $x \in M(k)$;
these can be also thought of as block matrices with $m$ identical $k \times
k$ blocks
along the diagonal.

Let us comment here that the homogeneous space obtained from (1) is a
``generalized Stieffel manifold" of ``orthogonal frames" of subspaces of
${\bf C}^n$ (or ${\bf R}^n$) with given pattern of dimensions.
 In both cases (1) and (2) the subgroup $H$ (resp. the Lie algebra
$\cal{H}$ consists of
(all) unitaries
(resp. skew-symmetric matrices having form (1) or (2)),  and the conditional
expectation is a norm one (with respect  to any unitarily invariant norm)
projection from $\cal{G}$,  the Lie algebra of $G = U(n)$, onto $\cal{H}$,  in
particular $\kappa(G/H) \leq 2$.  However,  except for the case $m = 2$ in
(1) (i.e. the Grassmann manifold) or (2),  we have $\kappa(G/H) > 1$.

Concerning the other parameters,  it is easily seen that in all cases
$\theta(G/H) = 2$ and $\pi/2 \leq diam \, G/H \leq \pi$.  Accordingly,
Theorem 8
gives good lower estimates for the covering numbers of $G/H$ and it remains to
handle the upper ones.  We will use the following (ad hoc)

\medskip\noindent {\bf Theorem 11.} {\em Let } $\alpha \in (0, 1/2]$;
{\em let } $n$, $G$ ($= U(n)$ or $SO(n)$),
$H$, $M = G/H$, $\cal{G}$, $\cal{H}$, {\em and } $d$ {\em be as before and
assume that }
$$\min \{ \theta (M), diam\, M, \kappa(M)^{-1} \} \geq \alpha.$$
{\em Furthermore,  asuume that one of the following holds }
\newline (a) $\dim\, H \leq (1 - \alpha) \dim\, G$
\newline (b) $H$ {\em acts reducibly on } ${\bf C}^n$ {\em (resp. } ${\bf
R}^n$
{\em and there is a reducing subspace } $E$ {\em with }
$\alpha n \leq \dim\, E \leq (1 - \alpha) n$
\newline (c) $H$ {\em acts reducibly on } ${\bf C}^n$ {\em (resp. } ${\bf
R}^n)$
{\em and there is a reducing subspace } $E$ {\em with }
$\dim\, E \equiv k \geq \alpha n$ {\em and such that the orthogonal
decomposition }
${\bf C}^n = E \oplus E^{\perp}$ {\em induces an isomorphism }
$H \rightarrow U(k) \times H_0$ {\em for some subgroup } $H_0$ {\em of }
$U(n-k)$ {\em (resp. } $ {\bf R}^n , SO(k), SO(n-k))$.
\newline {\em Then, for any } $\varepsilon \in (0,diam\, M]$,
$$(\frac {c}{\varepsilon} )^d \leq N(M, \varepsilon) \leq (\frac
{C}{\varepsilon} )^d,$$
\noindent {\em where } $c, C > 0$ {\em are constants depending only on }
$\alpha$.

\medskip\noindent {\bf Corollary 12.} {\em If } $H \subset U(n)$ {\em is
the group of
unitaries of a ``block-diagonal" or ``tensor-factor" algebra (described in
(1) or (2)),
then the assertion of Theorem 11 holds (e.g., with constants corresponding to }
$\alpha = 1/3$ {\em ).}

\medskip\noindent {\bf Proof.} We may of course assume $m \geq 2$. In the
case (2) the
condition (b) of Theorem 11 is always satisfied (with $\alpha = 1/3$). The
same is
true in the case (1) except if one of the projections $P_j$ is of rank
$\geq n/3$,
in which case (c) holds.

\medskip\noindent {\bf Remarks.} (i) As was pointed out in \cite{voiculescu},
Remark 7.2, the estimates for covering numbers  given by our Corollary 12
(the ``block-diagonal" case)
allow sharp free entropy and free entropy dimension estimates. Similarly,  the
 ``tensor-factor" case of the Corollary implies estimates for free entropy
and free entropy dimension of certain generators of free product von Neumann
algebras (\cite{voiculescupc}).
\newline (ii) The ``block-diagonal" case of Corollary 12 implies estimates for
covering numbers of some sets of matrices needed in \cite{dykema}.

\medskip\noindent {\bf Proof of Theorem 11.} As observed earlier,  it is
enough to
show the upper estimate.
\newline (a) The condition (a) is equivalent to $\dim\, M \geq \alpha
\dim\, G$.
It follows from Theorem 7 that, for any $\varepsilon \in (0, 2]$,  $G$ (hence,
by Lemma 2, $M = G/H$) admits an $\varepsilon$-net of cardinality
$ \leq (\frac {C}{\varepsilon} )^{dim \, G}$ .  If  $\varepsilon \geq
\beta$ and (a)
holds, this does not exceed $C(\alpha, \beta)^{\dim\, M}$  (the {\em real}
dimension).  It follows now from Lemma 10
that the image of $B_{\cal{X}}(r)$ (where $r = r(\alpha)$) contains a ball
in $M$
of radius $r_1 = r_1(\alpha)$, and so the former and the latter admit, for
$\varepsilon \leq r$,  an $\varepsilon$-net of cardinality
$ \leq (\frac {C_1 r}{\varepsilon} )^{dim \, M}$.  Combining this with the
preceding observation (applied to $\beta = r_1$) we get the required upper
estimate.
\newline (b) The condition (b)  implies (a)  (with $\alpha (1-\alpha)$ in place
of $\alpha$).
\newline (c) Since the arguments in the real and complex case are identical,
we restrict the discussion to the latter.
The condition (c) is not included in (b) only if
$\dim\, E \equiv k > (1 - \alpha) n$, in particular $k > n/2$.
Let $H_1 \subset G$ be the subgroup of the form $U(k) \times \{ I \}$
in the sense indicated in the condition (c),  then $H_1 \subset H$ and so
$M = G/H$ is a quotient of $M_1 = G/H_1$.  Now $M_1$ is isomorphic to
$G_{n,k} \times U(n-k)$ (for sure Lipschitz isomorphic with constant 2 if
the product metric on the latter is defined in the ``$\ell_{\infty}$ sense"),
and so, by Theorems 7 and 8, it admits,  for any $\varepsilon \in (0, 2]$,
an $\varepsilon$-net of cardinality
$\leq (\frac {C_2}{\varepsilon} )^{dim \, M_1}$
(again, the {\em real} dimension). Since $k > n/2$ implies
$dim \, M_1 < 2 dim \, M$,  arguing as in (a) we obtain the assertion.

In some applications (see e.g. \cite{exotic}, \cite{cube}) it is important
to know the metric entropy of $M$ equipped with a metric induced by
unitarily invariant norms other than the operator norm. The scheme presented
in this paper can be adapted to yield fairly sharp results in the general case.
Indeed, Lemmas 4, 5 and 6 involve statements about generic unitarily
invariant norms. Similarly, Lemma 10 and its proof carry over almost
word by word to the case of an arbitrary unitarily invariant norm
$\| \cdot \|$ once the parameters such as $\theta$ and $\kappa$ are
properly interpreted: the balls $B_{\cal{H}}(\cdot)$, $B_{\cal{X}}(\cdot)$
are to remain to be defined by the operator norm $\| \cdot \|_{\infty}$,
but $\theta (M)$ has to be the distance between $I$ and
$H \backslash \exp (B_{\cal{H}}(\pi))$ in the intrinsic
metric on $G$ induced by $\| \cdot \|$;  $\kappa(M)$ may be
calculated using $\| \cdot \|$ (which results in a quantity not larger
than the one given by the operator norm,  in particular $\kappa(M)=1$
if $\| \cdot \|$ is the Hilbert Schmidt norm). The ``linearization"
procedure can be then implemented and the problem is reduced to
estimating covering numbers of balls in $\cal{X}$ in the operator norm
with respect to $\| \cdot \|$.  As indicated in section 2 (see Lemma 1 and
the paragraph preceding it),  there exist numerous tools for obtaining
such estimates. In particular,  in many natural cases (e.g. $M = U(n)$,
$SO(n)$ or
the Grassmann manifolds $G_{n,k}$),  the volumes of a ball in $X$ with respect
to a unitarily invariant norm $\| \cdot \|$ and the inscribed operator norm
ball differ
by a factor $C^d$, $C$ - an universal constant  (this is easily implied
e.g., by classical
facts from \cite{Santalo}, \cite{chevet}, cf. \cite{fdbp}, p. 162;
see also \cite{gluskin} or \cite{SR}),  which allows to use Lemma 1
to show that $\log N(M,\varepsilon) \approx d \, \log(diam \, M
/\varepsilon)$  (and
$diam \, M \approx \| I \|$ if $M = U(n)$ or $SO(n)$,
$diam \, M \approx \| P \|$,  where $P$ is an orthogonal projection of rank
equal to $\min \{ k, n-k \}$ if $M = G_{n,k}$);  these cases have been
worked out in
\cite{Iowa}.


\begin{thebibliography}{00}

\bibitem{bhatia} R. Bhatia, \underline{Perturbation bounds for matrix
eigenvalues},
Pitman Research Notes 162,
Longman Scientific \& Technical, Harlow 1987.

\bibitem{busemann} H. Busemann,  \underline{Metric methods in Finsler
spaces and
the foundations of geometry}, Annals Math. Studies 8,
Princeton University Press, Princeton, 1942.

\bibitem{carl} B. Carl and I. Stephani, \underline{Entropy, Compactness and
Approximation of Operators}, Cambridge University Press, Cambridge,
1990.

\bibitem{chevet} S. Chevet, {\em S\'{e}ries de variables alleatoires
gaussiennes
\'{a} valeurs dans $E \hat{\otimes}_{\varepsilon} F$},
S\'{e}minaire sur la geometrie des espaces de Banach 1977-78, Ecole
Polytechnique, Palaiseau.

\bibitem{dykema} K. Dykema, {\em personal communication}.

\bibitem{gluskin} E. D. Gluskin, {\em The diameter of the Minkowski
compactum is roughly equal to } n,
Functional. Anal. Appl. 15 (1981),  72-73.

\bibitem{gohberg} I. C. Gohberg and M. G. Krein, \underline{Introduction to
the Theory of Linear}
\underline{Nonselfadjoint Operators}, Nauka, Moscow 1965.  English transl.:
AMS 1969

\bibitem{helgason} S. Helgason, \underline{Differential Geometry,
Lie Groups and Symmetric Spaces},
Academic Press, New York 1978.

\bibitem{herrero} D. Herrero and S. J. Szarek,  {\em How well can an
$n \times n$  matrix be approximated
by reducible ones?},  Duke Math. J. 53 (1986),  	233-248.

\bibitem{pajor} A. Pajor,  {\em personal communication}.

\bibitem{pisier} G. Pisier,
\underline{The Volume of Convex Bodies and Banach Space Geometry},
Cambridge University Press, Cambridge 1989.

\bibitem{SR} J. Saint-Raymond,  {\em Sur le volume des ideaux d'operateurs
classiques},
Studia Math. 80 (1984), 63-75.

\bibitem{Santalo} L.A. Santal\'{o}, {\em Un invariante afin para los cuerpos
convexos del espacio de } n {\em dimensiones}, Port. Math. 8(1949), 155-161.

\bibitem{Iowa} S. J. Szarek,  {\em Nets of Grassmann manifold and
orthogonal groups},  Proceedings of
Banach Space Workshop,  	University of Iowa Press
1982,  169-185.

\bibitem{fdbp} S. J. Szarek, {\em The finite dimensional basis problem with an
appendix on nets of Grassmann manifold},  Acta Math.  151 (1983),  153-179.

\bibitem{exotic} S. J. Szarek, {\em An exotic quasidiagonal operator},
J. Funct. Anal. 89 (1990), 274-290.

\bibitem{cube} S. J. Szarek, {\em Spaces with large distance to $
l_{\infty} ^n$ and random matrices},
Amer. J. Math. 112  (1990),  899-942.

\bibitem{tomczak} Nicole Tomczak-Jaegermann, \underline{Banach-Mazur
distances and finite-
dimensional}
 \underline{operator ideals}, Longman Scientific \& Technical,  Harlow 1989.

\bibitem{va} V. Varadarajan, \underline{Lie Groups, Lie Algebras and Their
Representations},
Springer Verlag, New York 1984.

\bibitem{voiculescu} D. Voiculescu,  {\em The analogues of entropy and of
Fisher's
information measure in free probability theory. III.  The absence of Cartan
subalgebras}, Geom. Funct. Anal. 6 (1996), 172-199.

\bibitem{voiculescupc} D. Voiculescu, {\em personal communication}.

\end{thebibliography}
\end{document}